\documentclass[psamsfonts, reqno]{amsart}

\usepackage{amscd,amsmath,amssymb,amsfonts,verbatim}
\usepackage{times}
\usepackage[cmtip, all]{xy}
\usepackage{graphicx}
\usepackage[active]{srcltx}
\usepackage{color}
\usepackage{textcomp}

\definecolor{refkey}{gray}{.5}   
\definecolor{labelkey}{gray}{.5} 
\definecolor{Red}{rgb}{1,0,0}

\newcommand{\pf}{{\bf Proof : }}

\newcommand{\qedwhite}{\hfill \ensuremath{\Box}} 

\newtheorem{theo}{Theorem}[section]
\newtheorem{prop}[theo]{Proposition}
\newtheorem{lem}[theo]{Lemma}
\newtheorem{cor}[theo]{Corollary}

\theoremstyle{definition}

\newtheorem{rem}[theo]{Remark}

\newtheorem{defi}[theo]{Definition}

\newcommand{\Um}{\mbox{\rm Um}}		\newcommand{\SL}{\mbox{\rm SL}}
\newcommand{\GL}{\mbox{\rm GL}}		
	\newcommand{\E}{\mbox{\rm E}}
\newcommand{\ESp}{\mbox{\rm ESp}}     \newcommand{\Sp}{\mbox{\rm Sp}}

\title{Nice group structure on the elementary orbit space of unimodular rows}
\author{Manoj K. Keshari and Sampat Sharma}
\newcommand{\Addresses}{{
  \bigskip
  \footnotesize

\textsc{Manoj K. Keshari, Department of Mathematics, IIT Bombay
            Mumbai 400076, INDIA}\par\nopagebreak
  \textit{E-mail:} Manoj K. Keshari \texttt{<keshari@math.iitb.ac.in>}

\medskip
  
 \textsc{Sampat Sharma, Department of Mathematics, IIT Bombay 
            Mumbai 400076, INDIA}\par\nopagebreak
  \textit{E-mail:} Sampat Sharma \texttt{<sampat@math.iitb.ac.in; sampat.iiserm@gmail.com>}

  \medskip

  }}

\begin{document}
\maketitle
\subjclass 2020 Mathematics Subject Classification:{13C10, 19D45, 19G12}

 ~ {Unimodular row, nice group structure}

 \begin{abstract}
\begin{enumerate}
\item Let $R$ be an affine algebra of dimension $d\geq 5$ over $\overline{\mathbb{F}}_{p}$ with $p>3$. Then the group structure on ${\Um_{d}(R)}/{\E_{d}(R)}$ is nice. 
\item Let $R$ be a commutative noetherian ring of dimension $d\geq 2$ such that $\mbox{E}_{d+1}(R)$ acts transitively on $\Um_{d+1}(R)$. Then the group structure on $\Um_{d+1}(R[X])/\mbox{E}_{d+1}(R[X])$ is nice.
 \end{enumerate}
 \end{abstract}
 
 \vskip0.50in

\section{Introduction} {\it All the rings are assumed to be commutative noetherian with unity $1$. }

 Let $n\geq 3$ be an integer and $R$ be a ring with stable dimension $\leq 2n-4$. It was proved by van der Kallen \cite{vdk2} that the orbit space ${\Um}_n(R)/{\E}_n(R)$ has  an abelian group structure. Here ${\Um}_{n}(R)$ and $\E_{n}(R)$ denote the set of unimodular rows and the group generated by elementary matrices (definition \ref{definitions}). If $v,w\in {\Um}_n(R)$, then in view of  {\cite[Lemma 3.2]{vdk3}},  the classes $[v]$ and $[w]$ can be expressed as $[v]=[(x_1,v_2,\ldots,v_n)]$ and $[w]=[(v_1,v_2,\ldots,v_n)]$. Then their product is defined in {\cite[Lemma 3.5 (i)]{vdk2}} as 
 $$[(x_1,v_2,\ldots,v_n)]\ast[(v_1,v_2,\ldots,v_n)]=[(v_1(x_1+w_1)-1,(x_1+w_1)v_2,v_3,\ldots,v_n)]$$
  where $w_1\in R$ is such that $v_1w_1=1$ modulo $(v_2,\ldots,v_n)$. We say that this group operation is {\bf{nice}} if it satisfies 
  $$ [(x_1,v_2,\ldots,v_n)]\ast[(v_1,v_2,\ldots,v_n)]=[(x_1v_1,v_2,\ldots,v_n)].$$ Let $C = {\mathbb{R}[x_{0}, x_{1},\ldots, x_{d}]}/{(x_{0}^{2} + x_{1}^{2} +\cdots + x_{d}^{2} - 1)}$ be the coordinate ring of the real $d$-sphere $S^{d}$ and $R= T^{-1}C,$ where $T$ is the multiplicative set of polynomial functions that do not have any zero on $S^{d}.$ Then the group structure on ${\Um}_{d+1}/{\E}_{d+1}(R)$ is not nice {\cite[Example 2.2(c)]{weibelci}}. We state some results where group structure on ${\Um}_n(R)/{\E}_n(R)$ is nice. 
\begin{enumerate}
\item If $R$ is a smooth affine algebra of dimension $d\geq 3$ over a perfect field $k$ with c.d.$_{2}(k)\leq 2$, then Fasel {\cite[Theorem 2.1]{faselbk}} proved that the group structure on 
${\Um}_{d+1}(R)/{\E}_{d+1}(R)$ is nice.
\item Garge and Rao {\cite[Theorem 3.9]{gr}} removed smoothness assumption in Fasel's result under the additional assumption that char$(k)\neq 2$ and c.d.$(k)\leq 1.$ 
\item If $R$ is a smooth affine algebra of dimension $d\geq 3$ over an algebraically closed field, then Gupta, Garge and Rao {\cite[Theorem 6.8]{ggr}} proved that  the group structure on ${\Um}_{d}(R)/{\E}_{d}(R)$  is nice. 
\item If $R$ is a local ring of dimension $d$ with $2R = R$, then Garge and Rao {\cite[Theorem 5.1]{gr}}  proved that 
the group structure on ${\Um}_{d+1}(R[X])/{\E}_{d+1}(R[X])$ is nice.
\end{enumerate}
In {\cite[Proposition 3.5]{raojose}}, Jose and Rao proved that if $R$ is a reduced affine algebra of dimension $d\geq 2$ over an algebraically closed field $k$, then the group structure on $\Um_{d+1}(R)/\E_{d+1}(R)$ is nice.  If $(c_{0}, c_{1}, \ldots, c_{d})\in \Um_{d+1}(R),$ then by Swan's version of Bertini theorem \cite{swanbertini}, we may assume that $R/(c_{1},\ldots, c_{d})$ is a finite product of $k.$ Thus $c_{0} = d_0^{2}~\mbox{mod}~(c_{1},\ldots, c_{d}).$ Let us see that this implies that  the group structure on  $\Um_{d+1}(R)/\E_{d+1}(R)$ is nice. Let $ v = (a_{0}, a_{1}, \ldots, a_{d}), w = (a_{0}', a_{1}, \ldots, a_{d})\in \Um_{d+1}(R)$. Then there exist $b_{0}, b_{0}'\in R$ such that 
$$[v] = [(b_{0}^{2}, a_{1}, \ldots, a_{d})]~\mbox{and}~ [w] = [(b_{0}'^{2}, a_{1}, \ldots, a_{d})].$$ Thus $a_{0}a_{0}' = b_{0}^{2}b_{0}'^{2}~\mbox{mod}~(a_{1}, a_{2}, 
\ldots, a_{d}).$ Therefore, 
\begin{align*} [v]\ast [w] &=  
  [(b_{0}^{2}b_{0}'^{2}, a_{1}, \ldots, a_{d})]  =  [(a_{0}a_{0}', a_{1}, \ldots, a_{d})]
\end{align*} by  {\cite[Lemma 3.5(v)]{vdk2}}. If $R$ is an affine algebra of dimension $d\geq 3$ over 
 $\overline{\mathbb{F}}_{p}$, then Jose and Rao {\cite[Proposition 3.5]{raojose}} states that the same proof will show that ${\Um}_{d}(R)/{\E}_{d}(R)$ has a nice group structure. We observe that if $v = (a_0, a_1,\ldots, a_{d-1})\in \Um_{d}(R)$, then $R/(a_{1},\ldots, a_{d-1})$ is one dimensional. So it is not a finite product of $\overline{\mathbb{F}}_{p}$ and hence we do not have $a_{0} = b_0^{2}~\mbox{mod}~(a_{1},\ldots, a_{d-1}).$ Therefore their proof does not work in this case.  We revisit the method of Garge and Rao \cite{gr} and prove the following results. 
\begin{theo}
 Let $R$ be an affine algebra of dimension $d\geq 5$ over $\overline{\mathbb{F}}_{p}$ with $p>3.$ Then the group structure on ${{\Um}_{d}(R)}/{{\E}_{d}(R)}$ is nice.
\end{theo}

\begin{theo}
 Let $R$ be a ring of dimension $d\geq 2$ such that ${\E}_{d+1}(R)$ acts transitively on $\Um_{d+1}(R).$ Then the group structure on $\Um_{d+1}(R[X])/{\E}_{d+1}(R[X])$ is nice. 
\end{theo}

The above result applies to following two cases where $\mbox{E}_{d+1}(R)$ acts transitively on $\Um_{d+1}(R)$: (i) Jacobson radical of $R$ has height $\geq 1$, (ii) $R$ is an affine algebra over $\overline{\mathbb{F}}_{p}$.
\par Let $R$ be a ring with $\mbox{sdim}(R)\leq 2n - 4$ and $n\geq 3.$ Then van der Kallen {\cite[Theorem 4.1]{vdk2}} has shown that the group ${{\Um}_{n}(R)}/{{\E}_{n}(R)}$ is the universal weak $n$-Mennicke symbol group $\mbox{WMS}_{n}(R).$ When the group strucutre on ${{\Um}_{n}(R)}/{{\E}_{n}(R)}$ is nice, the universal weak $n$-Mennicke symbol group $\mbox{WMS}_{n}(R)$ coincides with universal $n$-Mennicke 
symbol group $\mbox{MS}_{n}(R).$

\section{Preliminaries}
\par 
\begin{defi}
\label{definitions} Let $R$ be a ring and $I$ be an ideal of $R$.
\begin{itemize}
\item  A vector 
 $v = (a_{1},\ldots, a_{r})\in R^r$ is unimodular if there exist $(b_{1},\ldots, b_{r})\in R^r$ with 
$\langle v ,w\rangle = \Sigma_{i = 1}^{r} a_{i}b_{i} = 1$. The set 
of unimodular rows in $R^r$ is denoted as $\Um_r(R)$.   
Let $e_{i} = (0, \ldots, 1, \ldots, 0)\in \Um_{r}(R),$ where $1$ is at the $i^{th}$ place.
The set $\Um_{r}(R,I)$ consist of unimodular vectors $v\in\Um(R)$ which are congruent to $e_1$ modulo $I$. 
\item The group $\E_r(R)$ is generated by elementary matrices
$e_{ij}(\lambda) = \mbox{I}_{r} + \lambda {\E}_{ij} \in \GL_r(R)$, where $\lambda \in R$ and for $i\neq j$, 
${\E}_{ij} \in \mbox{M}_{r}(R)$ whose $(i,j)^{th}$ entry is $1$ and all other entries are zero. Here $\mbox{I}_{r}\in \GL_r(R)$ is identity matrix. The group
${\E}_{r}(R)$ acts on $\Um_r(R)$ by right multiplication, i.e.  $v\cdot\varepsilon := v\varepsilon.$
  We say $v \sim w$ if there exists an $\varepsilon \in E_{r}(R)$ such that $v = w\varepsilon.$ One 
can check that $\sim$ is an equivalence relation. Let ${\Um_{r}(R)}/{\E_{r}(R)}$ denote the set of orbits of this action. 
The equivalence class of $v\in \Um_r(R)$ is denoted by $[v].$

\item Let $\psi_{1} = \big( \begin{smallmatrix}
                 0 & 1\\
                 -1 & 0\\
                \end{smallmatrix}\big)$ and $\psi_{n} = \psi_{n-1} \perp \psi_{1}$ for $n\geq 2$. Let $\sigma_{n}$ be the permutation of $\{1, \ldots, 2n\}$ given by $\sigma_{n}(2i) = 2i-1$ and $\sigma_{n}(2i-1) = 2i$.
The symplectic group $\mbox{Sp}_{2m}(R) $ is the subgroup of $\GL_{2m}(R)$ consisting of all 
$\alpha \in {\GL}_{2m}(R)$ such that $ \alpha^{t}\psi_{m}\alpha = \psi_{m}$.
 
The elementary symplectic group $\mbox{ESp}_{2m}(R)$ is the subgroup of $\mbox{Sp}_{2m}(R)$ generated by elementary symplectic matrices  
$$
se_{ij}(z)=
\begin{cases}
I_{2m} + z\E_{ij}, &\textit{if}~ i = \sigma_{m}(j)\\
 I_{2m} + z\E_{ij} - (-1)^{i+j}z\E_{\sigma_{m}(j)\sigma_{m}(i)}, &\textit{if}~ i\neq \sigma_{m}(j)
\end{cases}
$$
for $ 1\leq i \neq j\leq 2m$ and $z\in R.$

\item Let $\mbox{GL}_{n}(R,I)$ denote the kernel of the canonical mapping 
$\mbox{GL}_{n}(R)\longrightarrow \mbox{GL}_{n}\left({R}/{I}\right)$ and  $\mbox{SL}_{n}(R,I)=\mbox{SL}_n(R)
\cap \mbox{GL}_{n}(R,I)$.

\item The group $\mbox{E}_{n}(I)$ is the subgroup of $\mbox{E}_{n}(R)$ generated
 by elements $e_{ij}(x)$ with $x\in I.$ The relative elementary group $\mbox{E}_{n}(R,I)$ is the normal closure of $\mbox{E}_{n}(I)$ in $\mbox{E}_{n}(R).$

\item The group $\mbox{ESp}_{2m}(I)$ is the subgroup of $\mbox{ESp}_{2m}(R)$ 
 generated by
 elements $se_{ij}(x)$ with $x\in I.$ The relative elementary symplectic group ${\ESp}_{2m}(R,I)$ is the normal closure of ${\ESp}_{2m}(I)$ in ${\ESp}_{2m}(R)$.

\item We say stable range condition $sr_{n}(R)$ holds for $R$ if 
 given $(a_{1},a_{2},\ldots,a_{n+1})\in \Um_{n+1}(R),$ there exist $c_{i}\in R$ such that 
 $(a_{1}+ c_{1}a_{n+1}, a_{2}+c_{2}a_{n+1},\ldots,a_{n}+c_{n}a_{n+1})\in \Um_{n}(R).$

\item The stable range
of $R,$ denoted by $sr(R),$ is the least integer $n$ such that $sr_{n}(R)$ holds. The stable dimension of $R$ is
$\mbox{sdim}(R) = sr(R) - 1.$
 
\item Let $(a,b)\in \Um_{2}(R)$ and $c,d\in R$ be such that $ad-bc =1.$ The Mennicke symbol $\mbox{ms}(a,b)$ is defined as follows {\cite[Chapter 2]{mg}} : 
$$\mbox{ms}(a,b) = \mbox{class of}~  \begin{bmatrix}
                 a & b & 0 \\
                 c & d & 0\\
0 & 0 & 1\\
                \end{bmatrix}\in \SL_{3}(R)/\E_{3}(R).$$
Observe that if $\mbox{ms}({a_1, b_1})  = \mbox{ms}({a_2, b_2}),$ then there exists $\sigma \in \SL_{2}(R)\cap \E_{3}(R)$ such that 
$({a_1, b_1}) \sigma = ({a_2, b_2}).$ Further $\mbox{ms}(aa', b) = \mbox{ms}(a,b)\mbox{ms}(a',b).$

\item The excision ring $R\oplus I$ has coordinatewise addition. The multiplication
 is given as 
 $(r,i).(s, j) = (rs, rj+si+ij)$, where $r,s \in R$ and $i,j\in I.$
 The multiplicative identity is $(1,0).$ 

\end{itemize}
\end{defi}

\par We note a result from {\cite[Proposition 3.1]{keshari}}.

\begin{lem}
\label{anjanlike}
 Let $A$ be a ring and $R$ be an affine $A$-algebra of dimension $d$. If $I$ is an ideal of $R$, then $R\oplus I$ is also an affine $A$-algebra of dimension $d.$
 \end{lem}         
 
 \begin{defi}
 A ring homomorphism $\phi: B\twoheadrightarrow D$ has a section if there exists a ring homomorphism 
  $\gamma : D \hookrightarrow B$ so that $\phi \circ \gamma $ is the identity on $D.$ In this case we call $D$ to be a retract 
  of $B.$
 \end{defi}
 
 The following lemma is due to Suslin and is an easy consequence of {\cite[Lemma 4.3, Chapter 3]{mg}}.

 \begin{lem}
  \label{k2}
  Let $B,D$ be rings and $\pi: B \twoheadrightarrow D$ has a section. If $J = \mbox{ker}(\pi),$ then 
  ${\E}_{n}(B,J) = {\E}_{n}(B) \cap {\SL}_{n}(B,J)$ for  $n\geq 3.$
 \end{lem}

 \begin{rem}
 \label{omega}
   Let $R$ be a ring and $I$ be an ideal of $R$. There is a natural homomorphism $\omega : R\oplus I \rightarrow R$ 
   given by $(x,i) \mapsto x +i \in R.$ The map $\gamma : R \rightarrow R\oplus I$ given by $x \mapsto (x, 0)$ is a section of $\omega.$ Further, $(R\oplus I)/(0\oplus I) \simeq R$ via the projection map 
$\pi : R\oplus I\longrightarrow R.$
   
   Let 
   $v = (1+i_{1}, i_{2}, \ldots, i_{n}) \in \Um_{n}(R,I)$, where $i_{j}$'s are in $I.$ Then
  $\tilde{v} = ((1,i_{1}), (0,i_{2}), \ldots, (0,i_{n}))\in \Um_{n}(R\oplus I, 0\oplus I)$. We call $\tilde v$ to be a lift of $v.$ Note 
  that $\omega$ sends $\tilde{v}$ to $v.$
 \end{rem}

We state two results from \cite[Corollary 17.3, Corollary 18.1, Theorem 18.2]{7}.

\begin{theo}
\label{transit} Let $R$ be an affine $C$-algebra of dimension $d\geq 2$, where $C$ is either a subfield $F$ of  $\overline{\mathbb{F}}_{p}$ or $C=\mathbb Z$. Then 
\begin{itemize}
\item If $d =2,$ then ${\E}_{3}(R)$ acts transitively on $\Um_{3}(R).$
\item If $d\geq 3,$ then $\mbox{sr}(R) \leq d.$ As a consequence, ${\E}_{d+1}(R)$ acts transitively on $\Um_{d+1}(R).$
\end{itemize}
\end{theo}

\begin{cor}
\label{transitivecor}  Let $R$ be an affine algebra of dimension $d\geq 2$ over a subfield $F$ of  $\overline{\mathbb{F}}_{p}$ and $I$ be an ideal of $R$. Then ${\E}_{d+1}(R, I)$ acts transitively on $\Um_{d+1}(R, I).$
\end{cor}
${\pf}$ Let $v\in \Um_{d+1}(R, I).$ In view of Lemma \ref{anjanlike}, $R\oplus I$ is an affine algebra  of dimension $d$ over $F.$ By theorem \ref{transit}, there exists $\varepsilon_{1} \in \mbox{E}_{d+1}(R\oplus I)$ such that 
$\tilde{v}\varepsilon_{1} = e_{1}.$ Going modulo $0\oplus I,$ we 
have $e_{1}\overline{\varepsilon_{1}} = e_{1}$ with $\overline{\varepsilon_{1}} \in 
  \mbox{E}_{d+1}(R).$ Now replacing $\varepsilon_{1}$ by $\varepsilon_{1}(\overset{-}\varepsilon_{1})^{-1}$ and
  using Lemma \ref{k2}, we may assume that  
  $\varepsilon_{1} \in \mbox{E}_{d+1}(R\oplus I, 0\oplus I)$ satisfying 
  $\tilde{v}\varepsilon_{1} = e_{1}.$ Now applying $\omega$, defined in Remark \ref{omega},  to last equation
  we get $v\varepsilon = e_{1}$ for some $\varepsilon \in \mbox{E}_{d+1}(R, I).$ 
     $~~~~~~~~~~~~~~~~~~~~~~~~~~~~~~~~~~~~~~~~~~~~~~~~~~~~~~~~~~~~~~~~~~~~~~~~~~~~~~~~~~~~~~~~~~~~~~~~~~~~~~~~~
           ~~~~~\qedwhite$   
\par Next we note a result of Vaserstein  {\cite[Lemma 5.5]{7}}.
\begin{lem}
\label{2.11}
  Let $R$ be a ring. Then for any $m\geq 1,$
 ${\E}_{2m}(R)e_{1} = ({\Sp}_{2m}(R) \cap {\E}_{2m}(R))e_{1}.$
\end{lem}
 
\begin{rem}
\label{2.12}
 It was observed in {\cite[Lemma 2.13]{pr}}  that Vaserstein's proof actually shows that  
 ${\E}_{2m}(R)e_{1} = {\ESp}_{2m}(R)e_{1}.$ The relative version was proved in {\cite[Theorem 5.5]{pr}} (when $2R = R$) and in {\cite[Theorem 3.9]{anjan}} in general. 
\end{rem}

\begin{theo}
\label{pratyusharao} Let $R$ be a ring and $I$ be an ideal of $R$. If $v\in \Um_{2n}(R, I),$ then $v {\E}_{2n}(R, I) = v {\ESp}_{2n}(R, I)$ for $n\geq 2.$
\end{theo}

\begin{cor}
\label{symplectictransitivecor} Let $R$ be an affine algebra of dimension $3$ over a subfield $F$ of  $\overline{\mathbb{F}}_{p}$ which $2R = R$ and $I$ be an ideal of $R$. Then ${\ESp}_{4}(R, I)$ acts transitively on $\Um_{4}(R, I).$
\end{cor}
${\pf}$  It follows from theorem \ref{pratyusharao} and corollary \ref{transitivecor}.
$~~~~~~~~~~~~~~~~~~~~~~~~~~~~~~~~~~~~~~~~~~~~~~~~~~~~~~~~~~~~~~~~~~~~~~~~~~~~~~~~~~~~~~~~~~~~~~~~~~~~~~~~~
           ~~~~~\qedwhite$

Vaserstein \cite{6} showed that the natural map $\phi_{n, n+1} : {\mbox{Sp}_{2n}(R)}/{\mbox{ESp}_{2n}(R)}{\longrightarrow}
 {\mbox{Sp}_{2n+2}(R)}/{\mbox{ESp}_{2n+2}(R)}$
is bijective for $n\geq \mbox{dim}(R)+2.$ The injective stability bound was improved by Rao and Basu \cite{br} and Basu, Chattopadhyay and Rao \cite{pr2} as follows.
\begin{theo}\label{basuraop} Let $R$ be a smooth affine algebra over a perfect $C_{1}$-field, then the natural map $$\phi_{n, n+1} : {{\Sp}_{2n}(R)}/{{\ESp}_{2n}(R)}\longrightarrow {{\Sp}_{2n+2}(R)}/{{\ESp}_{2n+2}(R)}$$ is injective for $n\geq \mbox{dim}(R)+1.$
\end{theo}

The next result is due to Suslin-Vaserstein {\cite[corollary 7.4]{7}}.
\begin{lem}
\label{vascor} Let $R$ be a commutative ring and $v, v' \in \Um_{3}(R).$ If $v\beta = v'$ for some  $\beta \in {\SL}_{3}(R) \cap {\E}_{4}(R),$ then there exists $\beta_{1} \in {\E}_{3}(R)$ such that $v\beta_{1} = v'.$
\end{lem}

 We note a result of Fasel {\cite[Lemma 3.3]{faselbk}}.
\begin{lem}
\label{fasellemma} Let $S$ be a smooth affine surface over an algebraically closed field of characteristic $\neq 2,3.$ Then ${\SL}_{2}(S) \cap {\E}(S) = {\SL}_{2}(S)\cap {\E}_{3}(S) = {\SL}_{2}(S)\cap 
{\ESp}_{4}(S) = {\SL}_{2}(S)\cap {\ESp}(S).$
\end{lem}
 
 The next result is due to Garge and Rao {\cite[Lemma 3.5]{gr}} for $n = d+1,$. The same proof works in our case.

\begin{lem}
\label{reduced} Let $R$ be a ring of stable dimension $\leq 2n-4$ and $n\geq 3.$ If the group structure on ${\Um}_n(R_{\mbox{red}})/{\E}_n(R_{\mbox{red}})$ is nice, then the group structure on ${\Um}_n(R)/{\E}_n(R)$ is nice. 
\end{lem}

Next we note a slightly modified version of pre-stabilization result of van der Kallen {\cite[Theorem 2.2]{vdk2}}. We get this version by applying van der Kallen's result  {\cite[Theorem 2.2]{vdk2}} to the matrix $g^{t}$ and then taking the transpose afterwards.

\begin{theo}
\label{prestab} Let $R$ be a ring of stable dimension $\leq 2n-3$ with $n\geq 3$ and $I$ be an ideal of $R.$ Let $i,j\geq 0$ and $g\in \GL_{n+i}(R)\cap \E_{n+i+j+1}(R, I).$ Then there exist matrices $u,v,w, M$ with entries in $I$ and 
$q$ with entries in $R$ such that 
$$ \begin{bmatrix}
                 I_{i+1} + uq & v \\
                 wq & I_{n-1}+M\\
                \end{bmatrix}\in \E_{n+i}(R,I) ~\textit{and}~ \begin{bmatrix}
                 I_{j+1} + qu & qv \\
                 w & I_{n-1}+M\\
                \end{bmatrix}\in g\E_{n+j}(R,I).$$

\end{theo}

We note a local-global principle for unimodular rows {\cite[Theorem 2.3]{topcoeffie}}. 
\begin{theo}
\label{raolg} Let $R$ be a ring, $v\in \Um_{n}(R[X])$ and $n\geq 3.$ Suppose that for all $\mathfrak{m}\in \mbox{Maxspec}(R),$ $v\equiv v(0)~(\mbox{mod}~\E_{n}(R_{\mathfrak{m}}[X]).$ Then $v\equiv v(0)~(\mbox{mod}~\E_{n}(R[X]).$
\end{theo}

The next result of van der Kallen and Rao {\cite[Theorem 1]{kallenrao}} is about the drop in injective stability for $\mbox{SK}_{1}$ in the case of non-singular affine algebra. 

\begin{theo}
\label{kr} Let $R$ be a non-singular affine algebra of dimension $d\geq 2$ over a perfect $C_{1}$-field $k.$ Then
$\SL_{d+1}(R) \cap \E_{d+2}(R) = \E_{d+1}(R)$
i.e. a stably elementary $\sigma \in \SL_{d+1}(R)$ belongs to $\E_{d+1}(R).$ Consequently, the natural map 
$\SL_{r}(R)/\E_{r}(R)\longrightarrow \mbox{SK}_{1}(R)$
is an isomorphism for $r\geq d+1.$
\end{theo}

\section{Main results}

We start with a lifting lemma for smooth affine algebras over a subfield $F$ of  $\overline{\mathbb{F}}_{p}$.

\begin{lem}
\label{lifting} Let $R$ be a three-dimensional smooth affine algebra over a subfield $F$ of  $\overline{\mathbb{F}}_{p}$ with $2R = R.$ Let $a\in R$ be such that $\mbox{dim}(R/aR) =2$ and $R/aR$ is smooth. If $\alpha \in {\SL}_{2}(R/aR) \cap {\ESp}(R/aR),$ then there exists 
$\beta \in {\SL}_{2}(R) \cap {\E}(R)$ such that $\alpha \equiv \beta ~\mbox{mod}~aR.$
\end{lem}
${\pf}$ By stabilization theorem \ref{basuraop}, $\alpha \in \mbox{ESp}_{6}(R/aR).$ Therefore there exists $\beta' \in \mbox{ESp}_{6}(R)$ such that 
$$\beta' ~\mbox{mod}~aR =  \begin{bmatrix}
                 I_{4} & 0\\
                 0 & \alpha \\
                \end{bmatrix}. $$
In view of Suslin-Vaserstein's estimates and theorem \ref{pratyusharao}, $\mbox{ESp}_{6}(R, aR)$ acts transitively on $\Um_{6}(R, aR).$ Therefore by altering $\beta'$ we may assume that its first row is $(1,0,0,0,0,0).$ By the symplectic condition, 
 $$\beta' =  \begin{bmatrix}
                 1 & 0 & 0 \\
                 \ast & 1 & \ast \\
 \ast & 0 & \beta'' \\
                \end{bmatrix}$$ for some $\beta'' \in \mbox{SL}_{4}(R).$
Now observe that $e_{1}\beta'' \in \Um_{4}(R, aR).$ In view of corollary \ref{symplectictransitivecor}, $\mbox{ESp}_{4}(R, aR)$ acts transitively on $\Um_{4}(R, aR).$ Therefore by altering $\beta'$ by an elementary symplectic matrix we may assume that
the first row of $\beta''$ is $(1,0,0,0).$ Now by the 
symplectic condition,
$$\beta' =  \begin{bmatrix}
                 1 & 0 & 0 & 0 & 0 \\
                 \ast & 1 & \ast & \ast & \ast \\
 \ast & 0 & 1 & 0 & 0 \\
\ast & 0 & \ast & 1 & \ast \\
\ast & 0 & \ast & 0 & \beta 
                \end{bmatrix}$$ 
and $\beta$ satisfies our requirements.
$~~~~~~~~~~~~~~~~~~~~~~~~~~~~~~~~~~~~~~~~~~~~~~~~~~~~~~~~~~~~~~~~~~~~~~~~~~~~~~~~~~~~~~~~~~~~~~~~~~~~~~~~~
           ~~~~~\qedwhite$

\begin{lem}
\label{mainlemma}  Let $R$ be an affine algebra of dimension $d\geq 4$ over a subfield $F$ of  $\overline{\mathbb{F}}_{p}$ with $2R = R$ and $(a_{1}, \ldots, a_{d})\in \Um_{d}(R)$. Let $B = R/(a_{4}R + \cdots + a_{d}R)$ and $C = R/(a_{3}R + \cdots + a_{d}R).$ Assume that $B$ and $C$ are smooth affine algebras over $F$ with $\mbox{dim}(B) \leq 3$ and $\mbox{dim}(C)\leq 2.$ Let $\alpha \in {\SL}_{2}(C)\cap {\ESp}(C)$ and 
$(\overline{a_{1}}, \overline{a_{2}})\cdot \alpha = (\overline{b_{1}}, \overline{b_{2}})$ for some $b_1,b_2\in R$. Then there exists $\gamma \in {\E}_{d}(R)$ such that $(a_{1}, \ldots, a_{d})\gamma = (b_{1}, b_{2}, a_{3}, \ldots, a_{d}).$
\end{lem}
${\pf}$ In view of lemma \ref{lifting}, there exists $\beta' \in \mbox{SL}_{2}(B)\cap \mbox{E}(B)$ such that $(\overline{a_{1}}, \overline{a_{2}}, \overline{a_{3}})\cdot ({\beta'\perp 1}) = (\overline{b_{1} + a_{3}u_{3}}, \overline{b_{2}+a_{3}v_{3}}, \overline{a_{3}})$ for some $\overline{u_{3}}, \overline{v_{3}}\in B.$ Using elementary transformations, we can find a $\beta \in \SL_{3}(B) \cap \E(B)$ such that $(\overline{a_{1}}, \overline{a_{2}}, \overline{a_{3}})\cdot {\beta} = (\overline{b_{1}}, \overline{b_{2}}, \overline{a_{3}}).$ By theorem \ref{kr}, $\beta \in \mbox{SL}_{3}(B)\cap \mbox{E}_{4}(B).$ Therefore by lemma \ref{vascor}, there exists $\beta_{1} \in \mbox{E}_{3}(B)$ such that 
$$(\overline{a_{1}}, \overline{a_{2}}, \overline{a_{3}})\cdot \beta_{1} = (\overline{b_{1}}, \overline{b_{2}}, \overline{a_{3}}).$$
Now our desired $\gamma$ can be easily constructed in the form 
$$\gamma =  \begin{bmatrix}
                 \beta_{1}' & 0 \\
                 \ast & I_{d-3} \\
 \end{bmatrix}$$ 
where $\beta_{1}'\in \mbox{E}_{3}(R)$ such that $\beta_{1}' \equiv \beta_{1} ~\mbox{mod}~(a_{4}R + \cdots + a_{d}R).$
$~~~~~~~~~~~~~~~~~~~~~~~~~~~~~~~~~~~~~~~~~~~~~~~~~~~~~~~~~~~~~~~~~~~~~~~~~~~~~~~~~~~~~~~~~~~~~~~~~~~~~~~~~
           ~~~~~\qedwhite$

The proof of next result is similar to {\cite[Proposition 3.8]{gr}}. We refer  to \cite{gr} for the definition of admissible pair and admissible sequence.
\begin{prop}
\label{mainprop} Let $R$ be an affine algebra of dimension $d\geq 5$ over a field $F\subset \overline{\mathbb{F}}_{p}$ with $p > 3.$ Let $I$ be an ideal of $R$ of height $\geq 1.$ Let $(v, w)$ be an admissible pair of size $d$ over $R.$ Then there exists another admissible pair 
$(v', w')$ of size $d$ over $R$ such that both $v', w' \equiv e_{d}~(\mbox{mod}~I)$ and $(v, w)\rightarrow (v', w')$ is an admissible sequence.
\end{prop}
${\pf}$ Let $v = (a, a_{1}, \ldots, a_{d-1})$ and $w = (b, a_{1}, \ldots, a_{d-1})$ be a pair of admissible rows of size $d$ over $R.$ Let bar denote reduction modulo $I$.  Let $\overline{v} = (\overline{a}, \overline{a_{1}}, \ldots, \overline{a_{d-1}})$ and $\overline{w} = (\overline{b}, \overline{a_{1}}, \ldots, \overline{a_{d-1}}).$ Since $\overline{v}$ is unimodular, there exist $\overline{f}, \overline{f_{1}},\ldots, \overline{f_{d-1}}\in \overline{R}$ such that $\overline{a}\overline{f} + \overline{a_{1}}\overline{f_{1}} + \cdots + \overline{a_{d-1}}\overline{f_{d-1}} = 1.$ In particular, the elements $\overline{a}$ and $ \overline{a_{1}}\overline{f_{1}} + \cdots + \overline{a_{d-1}}\overline{f_{d-1}}$ are coprime in $\overline{R}.$ Hence, by {\cite[Lemma 9.2]{7}}, there exists $\overline{z}\in \overline{R}$ such that $\overline{a'} = \overline{a} + \overline{a_{1}}\overline{z}\overline{f_{1}} + \cdots + \overline{a_{d-1}}\overline{z}\overline{f_{d-1}}$ is not a zero divisor. Thus $\mbox{ht}(\overline{a'})\geq 1.$ Similarly, one can find $\overline{y},  \overline{g_{1}},\ldots, \overline{g_{d-1}}\in \overline{R}$ such that if $\overline{b'} = \overline{b} + \overline{a_{1}}\overline{y}\overline{g_{1}} + \cdots + \overline{a_{d-1}}\overline{y}\overline{g_{d-1}},$ then  $\mbox{ht}(\overline{b'})\geq 1.$ Therefore  $\mbox{ht}(\overline{a'b'})\geq 1.$ Let $\overline{v_{1}} = (\overline{a'}, \overline{a_{1}}, \ldots, \overline{a_{d-1}})$ and $\overline{w_{1}} = (\overline{b'}, \overline{a_{1}}, \ldots, \overline{a_{d-1}}).$ 
Then $(\overline{v}, \overline{w}) \rightarrow (\overline{v_{1}} , \overline{w_{1}} )$ is an admissible sequence.

Let $\overline{\overline{R}}$ denotes reduction modulo $(\overline{a'b'})$ in $\overline{R}.$ Note that $(\overline{\overline{a_{1}}}, \ldots, \overline{\overline{a_{d-1}}})\in\Um_{d-1}(\overline{\overline{R}}).$ Since $\mbox{dim}(\overline{\overline{R}}) \leq d-2$ and $d-2\geq 3$, by theorem \ref{transit}, $\mbox{sr}(\overline{\overline{R}})\leq d-2.$ Therefore $\mbox{sr}_{d-2}(\overline{\overline{R}})$ holds. Thus there exist $\overline{\overline{c_{i}}} \in   \overline{\overline{R}}$ for $1\leq i\leq d-2$ such that  if 
$ \overline{\overline{a_{i}'}} = \overline{\overline{a_{i}}} + \overline{\overline{c_{i}}}\overline{\overline{a_{d-1}}}$ 
then $(\overline{\overline{a_{1}'}}, \ldots, \overline{\overline{a_{d-2}'}})\in \Um_{d-2}(\overline{\overline{R}}).$ Thus  $(\overline{a'}, \overline{a_{1}'}, \ldots, \overline{a_{d-2}'})\in \Um_{d-1}({\overline{R}})$ and $ (\overline{b'}, \overline{a_{1}'}, \ldots, \overline{a_{d-2}'})\in \Um_{d-1}({\overline{R}}).$ 

Summarizing the above discussion, we can modify the unimodular rows $\overline{v}$ and $\overline{w}$ to the unimodular rows $\overline{v_{2}} = (\overline{a'}, \overline{a_{1}'}, \ldots, \overline{a_{d-2}'}, \overline{a_{d-1}})$ and $ \overline{w_{2}} = (\overline{b'}, \overline{a_{1}'}, \ldots, \overline{a_{d-2}'}, \overline{a_{d-1}})$ such that $(\overline{a'}, \overline{a_{1}'}, \ldots, \overline{a_{d-2}'})\in \Um_{d-1}({\overline{R}})$ and $(\overline{b'}, \overline{a_{1}'}, \ldots, \overline{a_{d-2}'})\in \Um_{d-1}({\overline{R}}).$ Then $(\overline{v}, \overline{w}) \rightarrow (\overline{v_{1}} , \overline{w_{1}} ) \rightarrow (\overline{v_{2}} , \overline{w_{2}} )$ is an admissible sequence. Now follow the proof of {\cite[Proposition 3.8]{gr}} to complete the proof.

$~~~~~~~~~~~~~~~~~~~~~~~~~~~~~~~~~~~~~~~~~~~~~~~~~~~~~~~~~~~~~~~~~~~~~~~~~~~~~~~~~~~~~~~~~~~~~~~~~~~~~~~~~
           ~~~~~\qedwhite$

The proof of next result follows that of {\cite[Theorem 3.9]{gr}}. We will only mention the necessary changes.

\begin{theo}
\label{mainresult} Let $R$ be an affine algebra of dimension $d\geq 5$ over $\overline{\mathbb{F}}_{p}$ with $p > 3.$ Then the abelian group structure on the orbit space ${\Um_{d}(R)}/{\E_{d}(R)}$ is nice. 
\end{theo}
${\pf}$ In view of lemma \ref{reduced}, we may assume that $R$ is reduced. Let $v = (a,a_{1}, \ldots, a_{d-1}), w =(b,a_{1}, \ldots, a_{d-1}).$
\begin{enumerate}
\item Let $J$ be the ideal defining the singular locus of $R.$ Since $R$ is a reduced, $\mbox{ht}(J)\geq 1.$ By proposition \ref{mainprop}, we may assume 
$(v,w)$ to be an admissible pair such that $v, w \equiv e_{d}~(\mbox{mod}~J).$
\item Applying Swan's version of Bertini \cite{swanbertini} (see also \cite[P. 413]{murthy}), we can add a linear combination of $ab, a_{1}, \ldots, a_{d-2}$ to $a_{d-1}$ to get 
$a_{d-1}' = a_{d-1} + \lambda ab + \sum_{i = 1}^{d-2}\lambda_{i}a_{i}$ for $\lambda, \lambda_{i} \in R$ and ${R}/{(a_{d-1}')}$ is smooth outside the singular set of $R.$ Since $a_{d-1}' \equiv 1~(\mbox{mod}~J)$, ${R}/{(a_{d-1}')}$ 
is smooth of dimension atmost $d-1.$
\par Now with $b_{i} = a_{i} + \sum_{j =1}^{i-1}\lambda_{ji}a_{j} + \mu_{i}ab$ for suitable $\mu_{i}, \lambda_{ji}\in R$,  we have an admissible sequence $(v,w)\rightarrow (v', w')$ where 
$v' = (a, b_{1},\ldots, b_{d-2}, a_{d-1}'),$ $w' = (b, b_{1},\ldots, b_{d-2}, a_{d-1}')$ and $\mbox{ht}(b_{2}, \ldots, b_{d-2}, a_{d-1}')\geq d-2$, $\mbox{ht}(b_{3}, \ldots, b_{d-2}, a_{d-1}')\geq d-3.$ 
Note that $v', w' \equiv e_{d}~(\mbox{mod}~J).$ We also have $\mbox{dim}({R}/{(b_{2}, \ldots, b_{d-2}, a_{d-1}'))}\leq 2$ and 
$\mbox{dim}({R}/{(b_{3}, \ldots, b_{d-2}, a_{d-1}'))}\leq 3.$ By Swan's version of Bertini as in \cite{swanbertini},  ${R}/{(b_{2}, \ldots, b_{d-2}, a_{d-1}')}$ and ${R}/{(b_{3}, \ldots, b_{d-2}, a_{d-1}')}$ are smooth outside the singular set of $R.$ Since 
 $a_{d-1}' \equiv 1~(\mbox{mod}~J),$ both ${R}/{(b_{2}, \ldots, b_{d-2}, a_{d-1}')}$ and ${R}/{(b_{3}, \ldots, b_{d-2}, a_{d-1}')}$ are smooth. Note that $[v]=[v']$ and $[w]=[w']$
\item Now van der Kallen's product formula in the group ${\Um_{d}(R)}/{\E_{d}(R)}$ shows that 
$$ [v] \ast [w]=[v']\ast[w'] = [(a(b+p)-1, (b+p)b_{1},\ldots, b_{d-2}, a_{d-1}')]$$
where $p$ is choosen so that $ap-1$ belongs to the ideal generated by $b_{1}, \ldots, b_{d-2}, a_{d-1}'.$ Let `overline' denote the image in ${R}/{(b_{2}, \ldots, b_{d-2}, a_{d-1}')R}$. If `ms' denote the Mennicke symbol, then 
\begin{align*}
\mbox{ms}(\overline{a}(\overline{b}+\overline{p})-1, (\overline{b}+\overline{p})\overline{b}_{1}) & = \mbox{ms}(\overline{a}(\overline{b}+\overline{p})-1, (\overline{b}+\overline{p}))\mbox{ms}(\overline{a}(\overline{b}+\overline{p})-1, \overline{b}_{1}) \\
& = \mbox{ms}(\overline{a}(\overline{b}+\overline{p})-1, \overline{b}_{1}). 
\end{align*}
Therefore there is a $\overline{\sigma}\in \mbox{SL}_{2}(\overline{R})\cap \mbox{E}_{3}(\overline{R})$
 such that 
$$(\overline{a}(\overline{b}+\overline{p})-1, (\overline{b}+\overline{p})\overline{b}_{1})\overline{\sigma} = (\overline{a}(\overline{b}+\overline{p})-1, \overline{b}_{1}).$$ 
Therefore by lemma \ref{fasellemma}, $\overline{\sigma} \in { \mbox{ESp}(\overline{R})}.$ Now in view of lemma \ref{mainlemma}, there exists $\gamma \in E_{d}(R)$ such that 
$$({a}({b}+{p})-1, ({b}+{p}){b}_{1}, b_{2},\ldots, b_{d-2}, a_{d-1}'){\gamma} = ({a}({b}+{p})-1, {b}_{1}, b_{2},\ldots, b_{d-2}, a_{d-1}').$$ 
Since $ap-1$ belongs to the ideal generated by $b_{1}, \ldots, b_{d-2}, a_{d-1}',$ we have
\begin{align*}
[v'] \ast [w']  &= [({a}({b}+{p})-1, ({b}+{p}){b}_{1}, b_{2},\ldots, b_{d-2}, a_{d-1}')] \\
&=  [({a}({b}+{p})-1, {b}_{1}, b_{2},\ldots, b_{d-2}, a_{d-1}')]\\
&=   [({a}{b}, {b}_{1}, b_{2},\ldots, b_{d-2}, a_{d-1}')] \\
&= [({a}{b}, {a}_{1}, a_{2},\ldots, a_{d-2}, a_{d-1}')] \\
&= [({a}{b}, {a}_{1}, a_{2},\ldots, a_{d-2}, a_{d-1})]
\end{align*}
This proves the required product formula.

\end{enumerate}

$~~~~~~~~~~~~~~~~~~~~~~~~~~~~~~~~~~~~~~~~~~~~~~~~~~~~~~~~~~~~~~~~~~~~~~~~~~~~~~~~~~~~~~~~~~~~~~~~~~~~~~~~~
           ~~~~~\qedwhite$

 \begin{cor}
           \label{vandercor}
            Let $R$ be an affine algebra of dimension $d\geq 5$ over $\overline{\mathbb{F}}_{p}$ with $p > 3$. If $\sigma \in {\SL}_{d}(R)\cap {\E}_{d+1}(R)$, then 
          $e_{1}\sigma$ can be completed to an elementary matrix.
           \end{cor}
    ${\pf}$ In view of  theorem \ref{prestab}, we have
$$ [\sigma] = \begin{bmatrix}
                 1 + ux & uy \\
                 z^{t} & M\\
                \end{bmatrix},
                ~~ [\varepsilon] = \begin{bmatrix}
                 1 + ux & y \\
                 uz^{t} & M\\
                \end{bmatrix},$$ for some $u,x\in R, y,z\in M_{1,d-1}(R), M\in M_{d-1,d-1}(R), 
                \varepsilon \in E_{d}(R)$. Since the group structure on ${\Um_{d}(R)}/{\E_{d}(R)}$ is nice, we have 
$$[e_{1}\sigma] = [1+ux, uy] = [1+ux, u]\ast [1+ux,y] = [1+ux, y] =[e_{1}].$$
The last equality $ [1+ux, y] =[e_{1}]$ follows from the existance of $\varepsilon$.
$~~~~~~~~~~~~~~~~~~~~~~~~~~~~~~~~~~~~~~~~~~~~~~~~~~~~~~~~~~~~~~~~~~~~~~~~~~~~~~~~~~~~~~~~~~~~~~~~~~~~~~~~~
           ~~~~~\qedwhite$

\section{A nice group structure over polynomial rings}
In \cite{gr}, Garge and Rao proved that if $R$ is a local ring of dimension $d$ with $2R = R,$ then the group structure on $\Um_{d+1}(R[X])/\mbox{E}_{d+1}(R[X])$ is nice. In this section we study the group structure on 
$\Um_{d+1}(R[X])/\mbox{E}_{d+1}(R[X]),$ where $R$ is a ring of dimension $d\geq 2$ such that $\mbox{E}_{d+1}(R)$ acts transitively on $\Um_{d+1}(R).$ 

\begin{theo}
\label{generalnice} Let $R$ be a ring of dimension $d\geq 2$ with $2R = R$ such that ${\E}_{d+1}(R)$ acts transitively on $\Um_{d+1}(R).$ Then
\begin{enumerate}
\item the group structure on $\Um_{d+1}(R[X])/{\E}_{d+1}(R[X])$ is nice.
\item If $\sigma \in {\SL}_{d+1}(R[X])\cap {\E}_{d+2}(R[X]),$ then $[e_{1}\sigma] = [e_{1}].$
\end{enumerate}
\end{theo}
${\pf}$ (1)  Let $[u(X)] = [(ab, a_{1},\ldots, a_{d})] - [(a,a_{1}, \ldots, a_{d})]  - [(b, a_{1}, \ldots, a_{d})].$ For every $\mathfrak{m}\in \mbox{Maxspec}(R)$, by {\cite[Theorem 5.1]{gr}},
$[u(X)]_{\mathfrak{m}} = [e_{1}].$ Therefore by local-global principle \ref{raolg}, $[u(X)] = [u(0)].$ Note that $u(0)\in \Um_{d+1}(R),$ by hypothesis of the theorem $[u(0)] = [e_{1}].$ Thus $[u(X)] = [e_{1}],$ i.e. 
$$ [(a,a_{1}, \ldots, a_{d})] \ast [(b, a_{1}, \ldots, a_{d})] = [(ab, a_{1},\ldots, a_{d})].$$

(2)   Let $\mathfrak{m} \in \mbox{Maxspec}(R).$ In view of  {\cite[Corollary 5.17]{gr}}, $[e_{1}\sigma]_{\mathfrak{m}} = [e_{1}].$ By local-global principle \ref{raolg}, we have 
$[e_{1}\sigma] = [e_{1}\sigma(0)].$ Since $e_{1}\sigma(0)\in \Um_{d+1}(R),$ $[e_{1}\sigma(0)] = [e_{1}].$ Thus  $[e_{1}\sigma] = [e_{1}].$

$~~~~~~~~~~~~~~~~~~~~~~~~~~~~~~~~~~~~~~~~~~~~~~~~~~~~~~~~~~~~~~~~~~~~~~~~~~~~~~~~~~~~~~~~~~~~~~~~~~~~~~~~~
           ~~~~~\qedwhite$

\begin{lem}
\label{elementaryjacobson} Let $R$ be a ring of dimension $d\geq 2$ and height of the Jacobson radical $ J(R) \geq 1.$ Then ${\E}_{d+1}(R)$ acts transitively on $\Um_{d+1}(R).$
\end{lem} 
${\pf}$ Let $v\in \Um_{d+1}(R)$ and $\overline{R} = R/J(R).$ Since $\mbox{dim}(\overline{R})\leq d-1,$ by Bass \cite{bass}, there exists $\overline{\varepsilon} \in \mbox{E}_{d+1}(\overline{R})$ such that 
$\overline{v}\overline{\varepsilon} = \overline{e_{1}}.$ Let $\varepsilon\in \mbox{E}_{d+1}(R)$ be a lift of $\overline{\varepsilon}.$ Then one has $v\varepsilon = (1+a_{0}, a_{1}, \ldots, a_{d})$ for some $a_{i}\in J(R).$ Note that $1+a_{0}$ is a unit in $R$, 
thus $[v] = [e_{1}].$
$~~~~~~~~~~~~~~~~~~~~~~~~~~~~~~~~~~~~~~~~~~~~~~~~~~~~~~~~~~~~~~~~~~~~~~~~~~~~~~~~~~~~~~~~~~~~~~~~~~~~~~~~~
           ~~~~~\qedwhite$

Next we note a result of Rao from {\cite[Corollary 2.5]{invent}}.
\begin{lem}
\label{raocor} Let $R$ be a ring of dimension $d\geq 2$ with ${d!}R =  R.$ Then every $v(X)\in \Um_{d+1}(R[X])$ is extended from $R,$ i.e. $v(X){\underset{{\SL}_{d+1}(R[X])}{\sim}} v(0).$
\end{lem}

\begin{cor}
\label{raotype}  Let $R$ be a ring of dimension $d\geq 2$ with ${d!}R = R$ and $\mbox{ht}(J(R))\geq 1$. Then every $v(X)\in \Um_{d+1}(R[X])$ is completable.
\end{cor}
${\pf}$ In view of lemma \ref{raocor}, $v(X){\underset{\mbox{SL}_{d+1}(R[X])}{\sim}} v(0).$ Since $v(0)\in \Um_{d+1}(R)$, by lemma \ref{elementaryjacobson} $v(0)$ is elementarily completable. Thus $v(X)$ is completable.
$~~~~~~~~~~~~~~~~~~~~~~~~~~~~~~~~~~~~~~~~~~~~~~~~~~~~~~~~~~~~~~~~~~~~~~~~~~~~~~~~~~~~~~~~~~~~~~~~~~~~~~~~~
           ~~~~~\qedwhite$

\begin{cor}
\label{gargetype}  Let $R$ be a ring of dimension $d\geq 4$ with $2R=R$. Assume that either $R$ is an affine $C$-algebra with $C=\mathbb{ Z}$ or $C$ is a subfield $F$ of $\overline{\mathbb{F}}_{p}$  or $\mbox{ht}(J(R))\geq 1$.
Then 
\begin{enumerate}
\item The group structure on $\Um_{d+1}(R[X])/{\E}_{d+1}(R[X])$ is nice.
 \item Let $\sigma \in {\SL}_{d+1}(R[X])\cap {\E}_{d+2}(R[X]).$ Then $[e_{1}\sigma] = [e_{1}].$
 \end{enumerate}
\end{cor}
${\pf}$ (1) Follows from theorem \ref{generalnice}, lemma \ref{elementaryjacobson} and theorem \ref{transit}.
(2) Follows from theorem \ref{generalnice}.
$~~~~~~~~~~~~~~~~~~~~~~~~~~~~~~~~~~~~~~~~~~~~~~~~~~~~~~~~~~~~~~~~~~~~~~~~~~~~~~~~~~~~~~~~~~~~~~~~~~~~~~~~~
           ~~~~~\qedwhite$ 

 \medskip
\noindent
{\bf Acknowledgement:} We thank the referee for insisting to clarify the proof of Prop \ref{mainprop}, which
made some points clearer. We also thank the referee for going through
the manuscript with great care. A detailed list of suggestions by the referee improved
the exposition considerably.

\Addresses

\end{document}